\documentclass[11pt]{article}
\usepackage{amsmath,amssymb}
\usepackage{theorem}
\usepackage{epsfig}
\title{A Derivation of the Catalan Numbers 
       from a \\
       Bijection between Permutations and Labeled Trees}
\author{Bennet Vance \\ ({\tt bennet.vance@dartmouth.edu})}
\newcommand{\comment}[1]{}
{\theorembodyfont{\rmfamily}
 \newtheorem{exa}{Example}}

\comment{
\newcommand{\binom}[2]{(\mbox{\small$%
\begin{array}{@{}c@{}}
 #1 \\[-1.0ex] #2
\end{array}$})}
}

\newcommand{\bvia}{{\it a}}
\newcommand{\bvib}{{\it b}}
\newcommand{\bvic}{{\it c}}
\newcommand{\bvid}{{\it d}}
\newcommand{\bvbia}{\emph{\bfseries a}}
\newcommand{\bvbib}{\emph{\bfseries b}}
\newcommand{\bvbic}{\emph{\bfseries c}}
\newcommand{\bvbicp}{{\bf\textit{c$^{\prime}$}}}
\newcommand{\bvbid}{\emph{\bfseries d}}
\newcommand{\mname}[1]{\mbox{\it{#1\/}}}
\newcommand{\LEAF}[1]{\ell_{#1}}
\newcommand{\TREE}[3]{(#1 \odot_{#2} #3)}
\newcommand{\PFXTREE}[3]{\phi_{#2}(#1,#3)}
\newcommand{\bvd}{d}

\newcommand{\embrace}[3]{\raisebox{-0.2ex}[0ex][0ex]{%
    \raisebox{#2}[0ex][0ex]{\makebox[1em]{#1}}%
    \raisebox{#2}[0ex][0ex]{
                            $\left\lgroup\parbox[c]{0.5em}%
                              {\rule{0.0mm}{#3}}%
                             \right.$
                           }%
                          }}
\begin{document}
\maketitle
\section{Overview}
Let $C_{n}$ denote the $n$th Catalan number,
which represents (among other things)
the number of distinct binary trees
that have $n$ undistinguished nodes and $n+1$ undistinguished leaves.
(Here the term \emph{node} shall refer to \emph{nonterminal} nodes only.)
Now imagine that in such a tree,
one assigns $n$ distinct labels to the nodes,
and $n+1$ distinct labels to the leaves;
let $D_{n}$ denote the number of possible trees
with the nodes and leaves so labeled.
Because there are $n!$ ways to label the nodes,
and, independently,
$(n+1)!$ ways to label the leaves,
we see that $D_{n} = C_{n} \cdot n! (n+1)!$,
and hence that $C_{n} = D_{n} / (n! (n+1)!)$.

This note presents a bijection
between permutations of length $2n$
and binary trees having $n$ labeled nodes and $n+1$ labeled leaves.
From the existence of this bijection one may infer that
$D_{n} = (2n)!$.
The familiar formulas
for the Catalan numbers follow directly:
$C_{n} = (2n)! / (n! (n+1)!) = \binom{2n}{n}/(n+1)$.
More generally, the presented bijection maps
permutations of length $\bvd n$ to labeled $n$-node $\bvd$-ary trees,
and vice versa.
By an argument analogous to that for binary case,
one may conclude that the number of distinct \emph{unlabeled} $\bvd$-ary trees
with $n$ nodes and $m = (\bvd-1)n+1$ leaves is
$(\bvd n)! / (n! m!) = \binom{\bvd n}{n}/m$.

Section~\ref{secconcept} below illustrates
the bijection for the case of binary trees.
Section~\ref{secperm2tree} gives an algorithm
that maps permutations to labeled $\bvd$-ary trees,
and Section~\ref{sectree2perm} gives the inverse algorithm.
Section~\ref{secabstract} gives a more abstract characterization
of the bijection.
Section~\ref{secconclude} concludes.

\section{Tree Construction using Permutations}
\label{secconcept}

Let us begin by considering how one might construct a labeled binary tree
from tree fragments, using a permutation to guide the construction.

\begin{figure}[h]
\begin{center}
\begin{tabular}{@{}c@{\hspace{0.10em}}c@{\hspace{0.10em}}c@{}} 
\epsfig{file=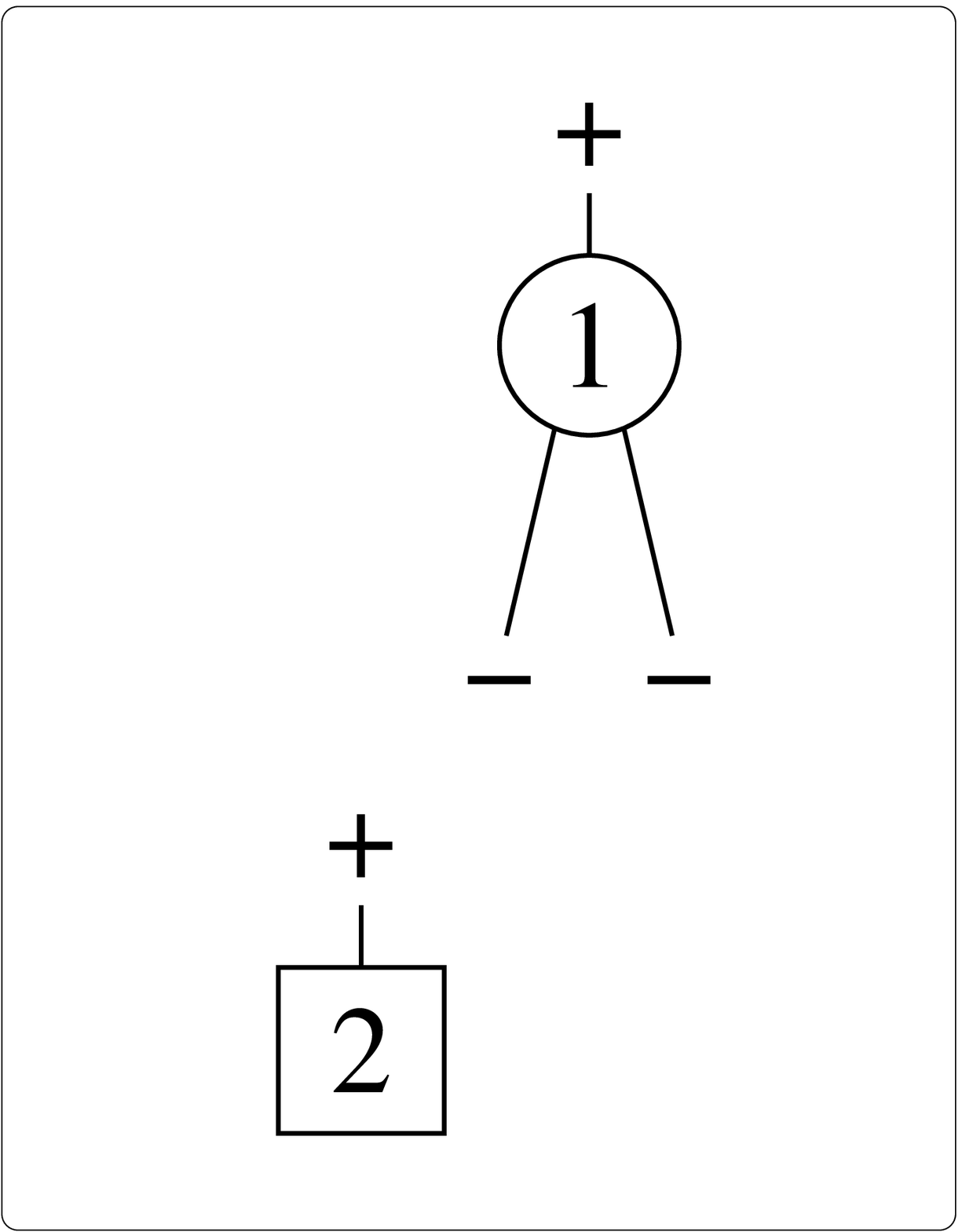,width=8.0em}&
\epsfig{file=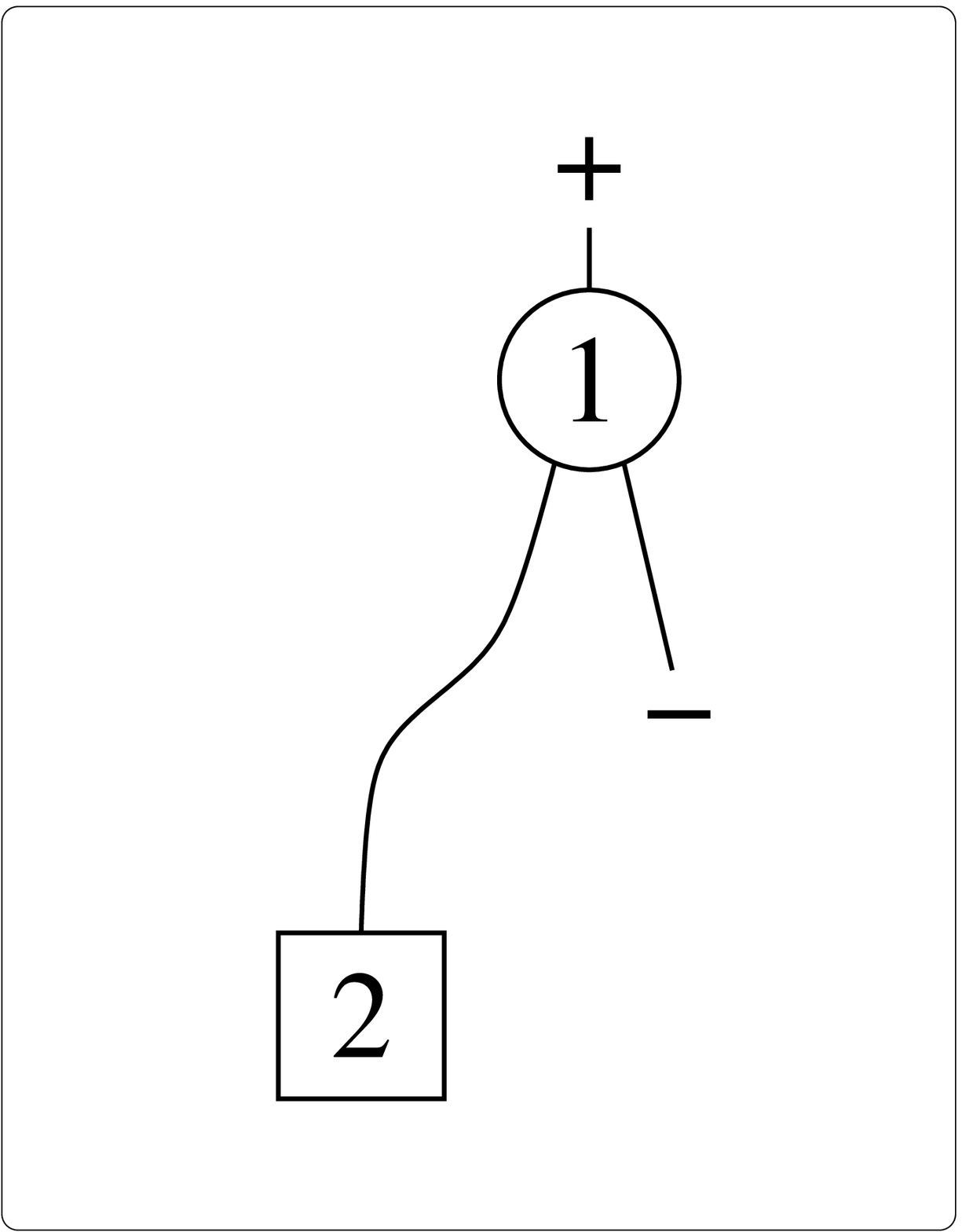,width=8.0em}&
\epsfig{file=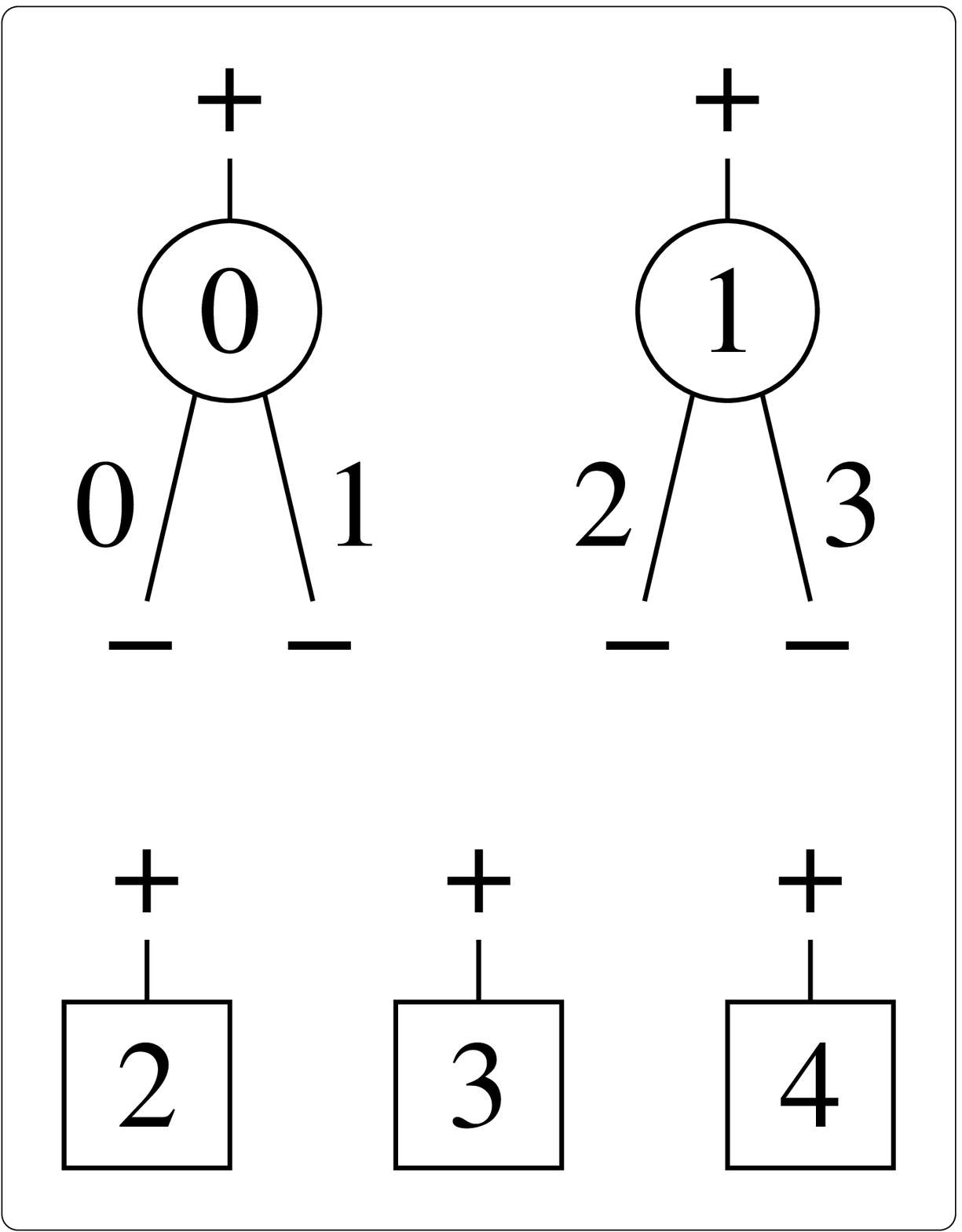,width=8.0em}\\
\bvia{} & \bvib{} & \bvic{}
\end{tabular}
\end{center}
\caption{Tree fragments}
\label{figfrags}
\end{figure}
Figure~\ref{figfrags}\bvia{}
depicts a leaf labeled $2$ and a node labeled $1$.
Each of these tree fragments has an \emph{upward sprout},
marked ``$+$'', directed toward an unresolved parent;
the node additionally has left and right \emph{downward} sprouts,
marked ``$-$'', directed toward unresolved children.
We associate with each sprout a unit \emph{charge} of the indicated sign.
The depicted leaf thus carries a net charge of $+1$, while
the node carries a net charge of $-1$
(i.e., one positive and two negative unit charges).
In the text, let $\LEAF{i}$ denote a leaf with label $i$,
and let $\TREE{u}{i}{v}$ denote a node with label $i$ and with
left and right subtrees $u$ and $v$, respectively;
let ``$-$'' denote an unresolved subtree.
Then Figure~\ref{figfrags}\bvia{} depicts $\LEAF{2}$ and $\TREE{-}{1}{-}$.

Figure~\ref{figfrags}\bvib{} shows the same leaf and node after
two of the sprouts have been connected by an edge.
The resultant structure is the incomplete tree $\TREE{\LEAF{2}}{1}{-}$.
When two sprouts of opposite sign are connected together,
their charges annihilate; total charge is thus conserved.
In the present instance, the leaf and node each acquire a net charge of zero
when they are connected.

In Figure~\ref{figfrags}\bvic{} we see the atomic fragments
needed to construct an arbitrary labeled $2$-node binary tree.
(We assume without loss of generality that a labeled $2$-node binary tree
has node labels $0$ and $1$, and leaf labels $2$, $3$, and $4$.)
The combined charge of all the fragments
comes to $+1$; 
moreover, since charge is conserved when fragments are connected,
any structure built from these fragments must also
carry a total charge of $+1$.
A consequence that we rely on below is that
any such structure
must contain at least one \emph{individual} node or leaf whose
net charge is positive, and hence $+1$
(since no other net positive charge can arise at a single node or leaf).
As illustrated in the figure,
to identify the left and right downward sprouts of node $i$
we use the numbers $2i$ and $2i+1$, respectively.
(We do not number \emph{upward} sprouts.)
Unlike node and leaf labels,
sprout numbers may not be freely reassigned:
they stand in a fixed relationship to the corresponding node labels.

\begin{figure}
\begin{center}
\begin{tabular}{@{}*{3}{c@{\hspace{0.10em}}}c@{}}
\epsfig{file=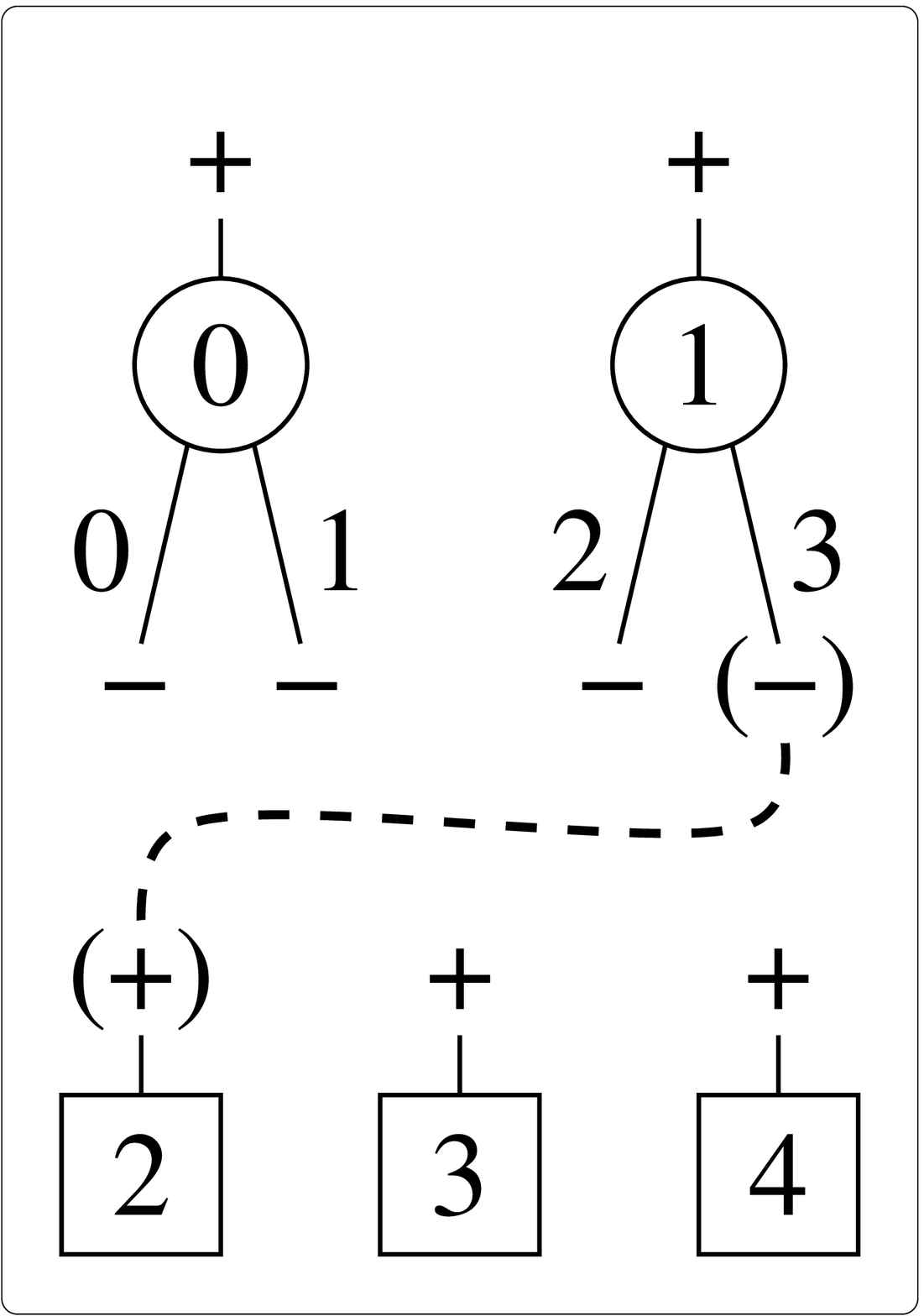,width=8.0em}&
\epsfig{file=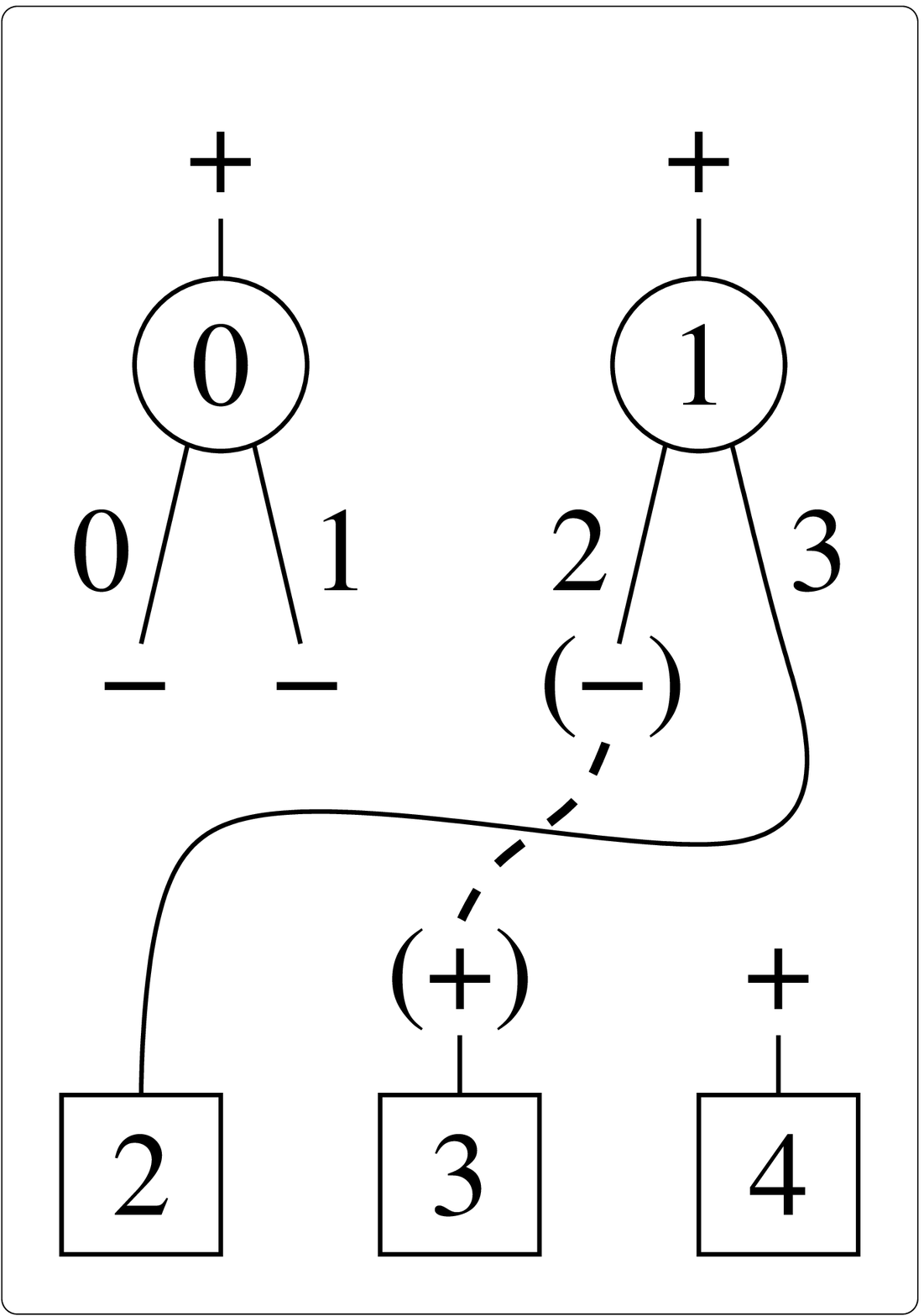,width=8.0em}&
\epsfig{file=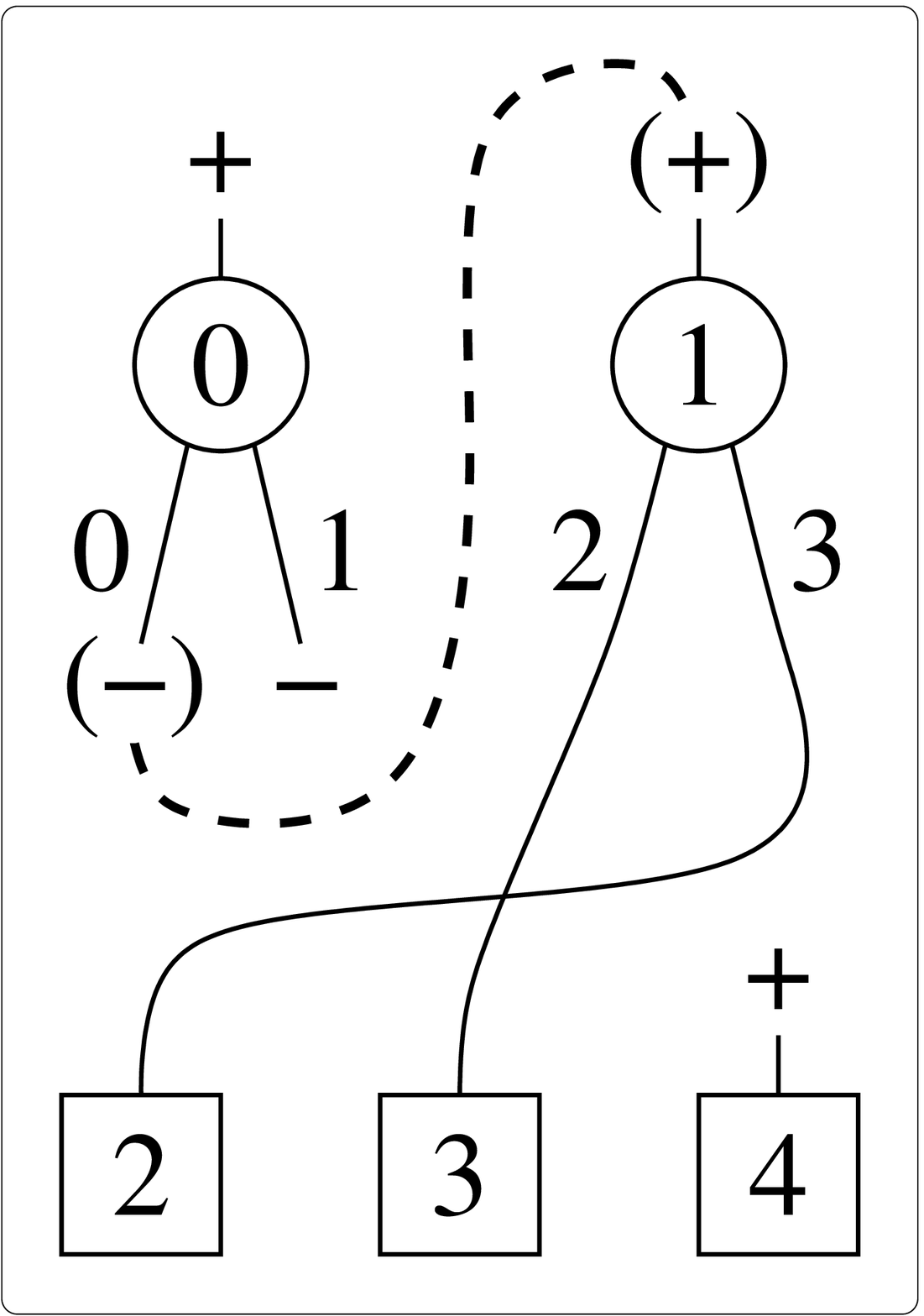,width=8.0em}&
\epsfig{file=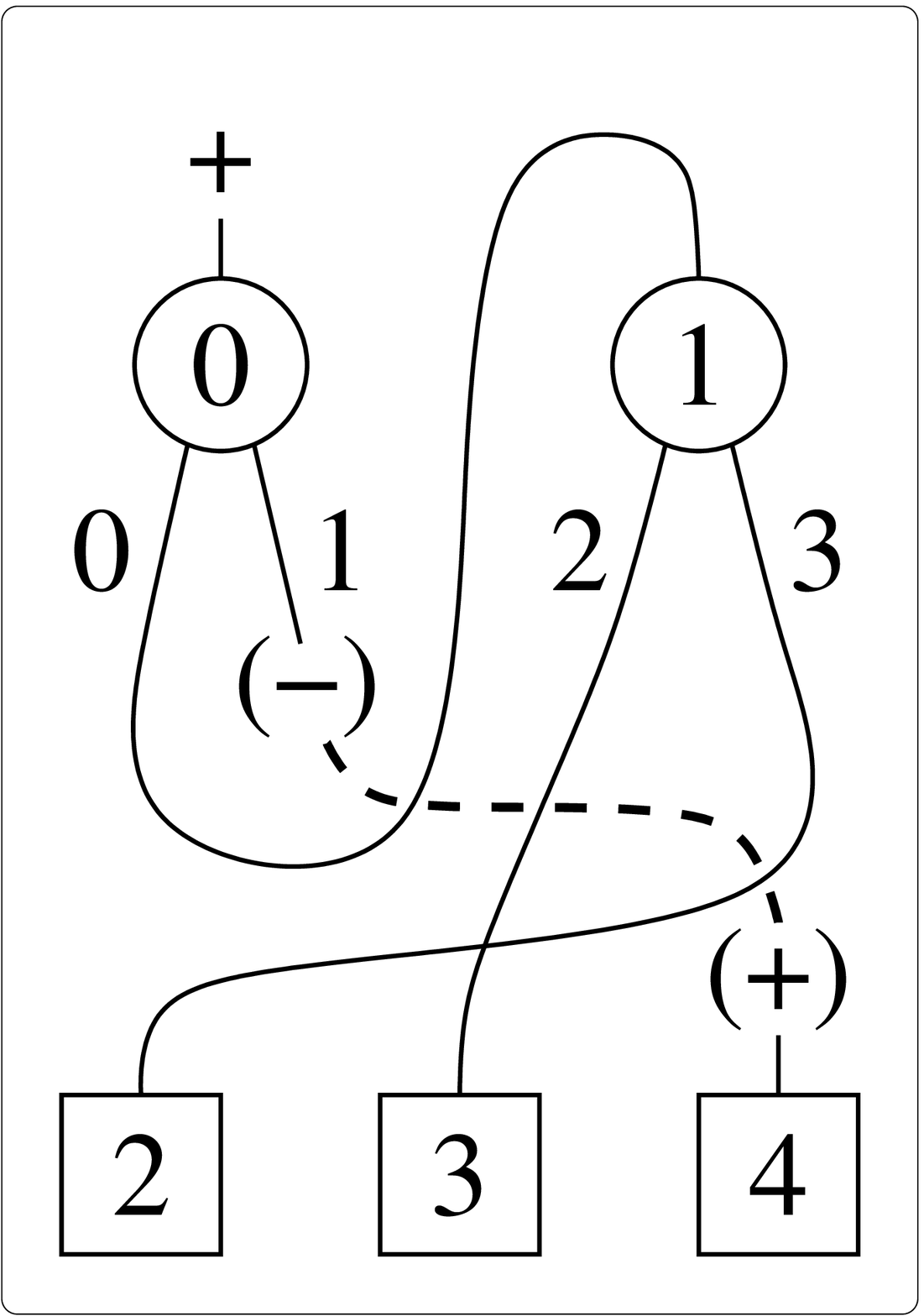,width=8.0em}\\
\bvia{} & \bvib{} & \bvic{} & \bvid{}
\end{tabular}
\end{center}
\caption{Step-by-step construction of a tree from fragments}
\label{figconstruct}
\end{figure}

Now suppose we wish to use a permutation---say,
$\langle 3, 2, 0, 1 \rangle$---to guide
the construction of a labeled $2$-node binary tree.
Figure~\ref{figconstruct} illustrates the step-by-step
construction of a tree from the fragments in Figure~\ref{figfrags}\bvic{}.
Each panel of Figure~\ref{figconstruct} augments the construction
with a new edge shown as a dashed line.
The downward sprouts participating in the new edges are numbered
$3$, $2$, $0$, and $1$, respectively;
in other words, we have used the permutation $\langle 3, 2, 0, 1 \rangle$
to decide the order in which the downward sprouts acquire children.
The child assigned to each new edge is determined as follows:
it is a node or leaf whose net charge is $+1$,
and among all such nodes and leaves, it is the one whose label is smallest.
(As noted above, at least one such node or leaf must always exist.)
Thus, in Figure~\ref{figconstruct}\bvia{} the new child is $\LEAF{2}$,
and in Figure~\ref{figconstruct}\bvib{} it is $\LEAF{3}$.
The node labels $0$ and $1$ are smaller than the leaf labels $2$ and $3$,
but it is only after acquiring the children $\LEAF{2}$ and $\LEAF{3}$
that one of the nodes---node~$1$---attains a net charge of $+1$.
Node~$1$ then becomes the new child in Figure~\ref{figconstruct}\bvic{};
$\LEAF{4}$ finally plays the role of new child
in Figure~\ref{figconstruct}\bvid{}.
At the end of this construction we have the tree
$\TREE{\TREE{\LEAF{3}}{1}{\LEAF{2}}}{0}{\LEAF{4}}$.
The completed tree can be characterized as a finite function
from downward sprout numbers to child labels,
namely $\{ 0 \mapsto 1, 1 \mapsto 4, 2 \mapsto 3, 3 \mapsto 2 \}$,
or, equivalently, as a vector of \emph{child pointers}
$\langle 1, 4, 3, 2 \rangle$.
This vector contains every node and leaf label except that of the tree's root
(node $0$ in this example).

Distinct permutations yield distinct labeled trees under this construction,
as we shall establish through the algorithms presented below.
Here we merely show that the illustrated construction must yield
a \emph{tree}, as opposed to a cyclic or disconnected graph of some kind.

Call a tree or subtree \emph{complete} if it has no unresolved descendants;
call a \emph{node} complete if it is the root of a complete tree or subtree.
Leaves are trivially complete.
Now imagine that the acquisition of a child $c$ by a parent $p$
results in a cycle among the tree fragments;
this outcome is possible only if $c$ was previously an ancestor of $p$.
But since $p$ was previously incomplete (else it could not have acquired the
child $c$), its ancestor $c$ was also incomplete.
Thus, a cycle can arise only when a node acquires an incomplete child.

At each step of the construction above we required that the new child 
be a node or leaf with net charge $+1$.
That requirement implied inductively
that each such node or leaf was complete:
it had no unconnected downward sprouts
(which would have contributed negative charge),
and its children, if any, must have themselves been complete
when it acquired them.
By ensuring completeness of new children,
the requirement of net charge $+1$ thus precluded cycles.
In addition, the presence of positive charge
indicated that a prospective child did not already have a parent,
so no separate bookkeeping was needed to make that determination.

Each time an edge was added to the construction,
two previously disconnected tree fragments became connected,
and the total number of disconnected fragments decreased by one.
(The two fragments joined by the edge could not have been connected previously,
because if they had been,
then connecting them anew would have created a cycle---a possibility
we have ruled out.)
After the addition of four edges, therefore, the five initial fragments
of Figure~\ref{figfrags}\bvic{} necessarily coalesced
into a connected acyclic graph.
The unique node carrying net charge $+1$ at the end of the construction
had to be the root of a complete tree which, by connectivity,
could be counted on to contain every node and leaf in the graph.

\section{A General Algorithm}
\label{secperm2tree}

The construction from Section~\ref{secconcept} above
extends straightforwardly to $\bvd$-ary trees.
Consider a complete $\bvd$-ary tree with $n$ nodes and $m$ leaves.
Each node is parent to $d$ edges, giving $\bvd n$ edges altogether.
In a construction that adds one edge connecting two fragments on each step,
the number of construction steps must equal the number of edges $\bvd n$,
and the initial number of fragments $n+m$ must exceed this value by $1$;
that is, $n+m = \bvd n + 1$.
This constraint on $\bvd$, $n$, and $m$ can be seen
in terms of charge as well:
Each of the $\bvd n$ construction steps consumes a positive charge
at one of the $n + m$ initial upward sprouts,
and a negative charge at one of the $\bvd n$ initial downward sprouts;
when the tree is fully connected,
only a single positive charge survives at the root.
We thus confirm the constraint
$(n + m - \bvd n) - (\bvd n - \bvd n) = +1$, or $n + m - \bvd n = +1$,
with the corollary $m = (\bvd - 1)n + 1$.

Figure~\ref{figperm2tree} gives a True BASIC program
that maps permutations of length $\bvd n$ to labeled $d$-ary trees.
\begin{figure}
\verb?    for i = 0 to d*n? \\
\verb?       if i <  n then let charge[i] = +1 + d*(-1)? \\
\verb?       if i >= n then let charge[i] = +1? \\
\verb?    next i? \\
\embrace{\bvbia}{7.45ex}{10.3ex} \\[-1.0ex]
\verb?    for i = 0 to d*n - 1? \\
\raisebox{0ex}[0ex][0ex]{} \\[-1.0ex]
\verb?       for j = 0 to d*n? \\
\verb?          if charge[j] = +1 then exit for? \\
\verb?       next j? \\
\embrace{\bvbib}{6.0ex}{7.4ex} \\[-1.0ex]
\verb?       let k = perm[i]? \\
\verb?       let kid[k] = j? \\
\embrace{\bvbic}{4.55ex}{4.5ex} \\[-1.0ex]
\verb?       let charge[j] = 0? \\
\verb?       let charge[int(k/d)] = charge[int(k/d)] + 1? \\
\embrace{\bvbid}{4.55ex}{4.5ex} \\[-1.0ex]
\verb?    next i?
\caption{Code for mapping a permutation to a labeled tree}
\label{figperm2tree}
\end{figure}
Assumed given are $\bvd$, $n$, and a permutation
$\mname{perm}[0]$, \ldots, $\mname{perm}[\bvd n - 1]$
of the integers $0, \ldots, \bvd n - 1$.
From these inputs the program computes
a representation of a labeled $n$-node $\bvd$-ary tree
as a vector of child pointers
$\mname{kid}[0]$, \ldots, $\mname{kid}[\bvd n - 1]$,
where $\mname{kid}[\bvd q + r]$ gives
the label of the $r$th child of node $q$,
for $0 \leq q < n$ and $0 \leq r < \bvd$.
The integers $0, \ldots, n-1$ serve as node labels,
and $n, \ldots, \bvd n$ as leaf labels.

The code segment marked {\bvbia} in Figure~\ref{figperm2tree}
initializes the array $\mname{charge}[0]$, \ldots, $\mname{charge}[\bvd n]$,
which gives, for each $i$,
the net charge of the node or leaf labeled $i$.
Each node starts out with net charge $1 - \bvd$
(reflecting one upward and $\bvd$ downward sprouts),
and each leaf, with $+1$ (reflecting just an upward sprout).
As before, we have the invariant
$\sum_{i} \mname{charge}[i] \equiv +1$,
which again assures us that at every step of our construction
there will be some node or leaf whose net charge is positive.

The remainder of Figure~\ref{figperm2tree} consists of a loop on $i$
that iteratively performs one construction step, as follows:
First (in code segment {\bvbib}),
it finds the smallest label $j$ such that $\mname{charge}[j] = +1$;
second (in code segment {\bvbic}),
it makes node or leaf $j$ the child of downward sprout number $k$,
where $k = \mname{perm}[i]$;
and third
(in code segment {\bvbid}),
it updates the charges associated with node or leaf $j$
and with its new parent.
The expression $\mbox{\tt int}(k/\bvd)$ in code segment {\bvbid}
denotes $\lfloor k/\bvd \rfloor$,
which yields the label of the parent node
that owns downward sprout number $k$.
When the code in Figure~\ref{figperm2tree} completes
after $\bvd n$ loop iterations,
the resultant tree's root label can be ascertained
by once again computing $\lfloor k/\bvd \rfloor$:
in general, this expression identifies a node that
has just acquired a child,
and the node that acquires a child last must be the root.

Figure~\ref{figrun24} shows
\begin{figure}
\[
\begin{array}{%
@{\mbox{\hspace{0em}}\langle}c@{,}c@{,}c@{,}c@{\rangle\;\leftrightarrow\;\;}c%
@{\mbox{\hspace{5em}}\langle}c@{,}c@{,}c@{,}c@{\rangle\;\leftrightarrow\;\;}c%
}
0& 1& 2& 3&  ((a + b) \times c) &     2& 0& 1& 3&  (a \times (b + c)) \\
0& 1& 3& 2&  (c \times (a + b)) &     2& 0& 3& 1&  (b + (a \times c)) \\
0& 2& 1& 3&  (b \times (a + c)) &     2& 1& 0& 3&  (a \times (c + b)) \\
0& 2& 3& 1&  (a + (b \times c)) &     2& 1& 3& 0&  ((a \times c) + b) \\
0& 3& 1& 2&  ((a + c) \times b) &     2& 3& 0& 1&  ((a \times b) + c) \\
0& 3& 2& 1&  (a + (c \times b)) &     2& 3& 1& 0&  (c + (a \times b)) \\[1ex]
1& 0& 2& 3&  ((b + a) \times c) &     3& 0& 1& 2&  ((b + c) \times a) \\
1& 0& 3& 2&  (c \times (b + a)) &     3& 0& 2& 1&  (b + (c \times a)) \\
1& 2& 0& 3&  (b \times (c + a)) &     3& 1& 0& 2&  ((c + b) \times a) \\
1& 2& 3& 0&  ((b \times c) + a) &     3& 1& 2& 0&  ((c \times a) + b) \\
1& 3& 0& 2&  ((c + a) \times b) &     3& 2& 0& 1&  ((b \times a) + c) \\
1& 3& 2& 0&  ((c \times b) + a) &     3& 2& 1& 0&  (c + (b \times a))
\end{array}
\]
\caption{A correspondence between permutations and labeled binary trees}
\label{figrun24}
\end{figure}
a correspondence between
permutations of length $4$ and labeled $2$-node binary trees.
This correspondence was obtained by running
the algorithm of Figure~\ref{figperm2tree} on each possible permutation
of $\langle 0,1,2,3 \rangle$ with $\bvd = 2$ and $n = 2$.
For readability, the nodes in Figure~\ref{figrun24}
are designated $+$ and $\times$, and the leaves, $a$, $b$, and $c$.
Observe that $a$ is the left-hand operand of $+$
in the first six expressions,
which correspond to permutations that start with $0$;
in the next six, $a$ is the \emph{right}-hand operand of $+$,
reflecting permutations that start with $1$;
and so on.
Similar patterns \emph{within} each group of six expressions
reflect the second, third, and fourth components of the permutations.

\section{Inverting the Algorithm}
\label{sectree2perm}

The algorithm of Figure~\ref{figperm2tree}
can easily be run ``in reverse'' to perform the inverse
mapping from labeled $\bvd$-ary trees to permutations.
That is, given $\bvd$, $n$,
and a representation of a tree in the $\mname{kid}[\,]$ array,
the algorithm can be made to retrace the construction of that tree, and
to compute in $\mname{perm}[\,]$
the permutation that would have caused
that particular tree to be constructed.

To achieve this reversal, one need only replace
code segment {\bvbic}  in   Figure~\ref{figperm2tree} with
code segment {\bvbicp} from Figure~\ref{figtree2perm}.
\begin{figure}[b]
\mbox{\qquad\qquad\qquad}\vdots \\
\raisebox{0ex}[0ex][0ex]{} \\[-1.0ex]
\verb?       for k = 0 to d*n - 1? \\
\verb?          if kid[k] = j then exit for? \\
\verb?       next k? \\
\verb?       let perm[i] = k? \\
\embrace{\bvbicp}{7.45ex}{10.3ex} \\[-1.0ex]
\mbox{\qquad\qquad\qquad}\vdots
\caption{Code revision for mapping a labeled tree to a permutation}
\label{figtree2perm}
\end{figure}
The difference between {\bvbic} and {\bvbicp}
is that where the former consults $\mname{perm}[\,]$
and assigns to $\mname{kid}[\,]$, the latter does the opposite:
it consults $\mname{kid}[\,]$ and assigns to $\mname{perm}[\,]$.
Neither array is referenced anywhere else in the algorithm,
so it is immaterial which array is the source of information,
and which the recipient.
The only aspect of code segment {\bvbic} or {\bvbicp}
that matters to the rest of the algorithm is that
(aside from not disturbing $i$, $j$, or $\mname{charge}[\,]$)
this code segment must furnish in $k$ a succession of distinct values
from $0, \ldots, \bvd n - 1$.
In the case of code segment {\bvbic},
this requirement is met in that
the values $k$ come from the permutation $\mname{perm}[\,]$.
In the case of {\bvbicp},
the values $k$ are the indices of the child pointers that
point to the labels $j$
of all the non-root nodes and leaves of the given tree;
since each such node or leaf is the target of a unique child pointer,
each value $k$ in $0, \ldots, \bvd n - 1$ will be furnished exactly once.
In either case, we have $\mname{perm}[i] = k$ and $\mname{kid}[k] = j$
after execution of code segment {\bvbic} or {\bvbicp}.
Thus, the relationship between
$\mname{perm}[\,]$ and $\mname{kid}[\,]$ will be the same
whether the algorithm is run forward or ``in reverse.''

Code segment {\bvbib}
in Figure~\ref{figperm2tree}
chose, from among the nodes and leaves carrying charge $+1$,
the one with the \emph{smallest} label $j$;
but that choice was somewhat arbitrary.
So long as code segment {\bvbib} implements \emph{some}
deterministic policy for
choosing a $j$ such that $\mname{charge}[j] = +1$,
the algorithm of Figure~\ref{figperm2tree}
remains invertible through substitution of
code segment
{\bvbicp} for {\bvbic}.

\section{An Abstract Characterization}
\label{secabstract}

The mappings presented above
also admit a more abstract characterization.
Assume $\bvd$ and $n$ fixed,
and let $Q = \{ 0, \ldots, n-1 \}$
and $R = \{ 0, \ldots, \bvd-1 \}$.
For $q$ in $Q$,
let $\phi_{q}(\ldots, u_{r},\ldots)$ denote a node labeled $q$
with immediate subtrees $u_{r}$ for $r$  in $R$.
If $d=2$, the notation
$\PFXTREE{u_{0}}{q}{u_{1}}$ is equivalent to $\TREE{u_{0}}{q}{u_{1}}$.
We define the $\mname{maxleaf}$ and $\mname{height}$ of a tree
in the obvious way:
\[
\begin{array}{l@{\qquad\qquad}l}
\mname{maxleaf}\;\LEAF{p} = p & \mname{height}\;\LEAF{p} = 0 \\
\mname{maxleaf}(\phi_{q}(\ldots, u_{r}, \ldots)) = &
\mname{height}(\phi_{q}(\ldots, u_{r}, \ldots)) = \\
\qquad         \max \{ \ldots, \mname{maxleaf}\;u_{r}, \ldots \} &
\qquad         1 + \max \{ \ldots, \mname{height}\;u_{r}, \ldots \}
\end{array}
\]
We then define the relation $u \prec v$ to hold on trees $u$ and $v$
just if $\mname{maxleaf}\;u < \mname{maxleaf}\;v$,
or if $\mname{maxleaf}\;u = \mname{maxleaf}\;v$
and $\mname{height}\;u < \mname{height}\;v$.

Now let $t$ be a $d$-ary tree
with node labels $0, \ldots, n-1$,
and leaf labels $n, \ldots, \bvd n$.
We define the permutation $P(t)$ as follows.
First,
for $k$ in $0, \ldots, \bvd n - 1$,
we define $\sigma_{k}(t)$,
or $\sigma_{k}$ for short,
to be the proper subtree $v$ of $t$ such that
$v$ is the $r$th child
of the node labeled $q$,
where $q = \lfloor k/\bvd \rfloor$ and $r = k - \bvd q$.
We then take $P(t)$ to be the permutation $\pi$ on $0, \ldots, \bvd n - 1$
such that $\sigma_{\pi(0)} \prec \ldots \prec \sigma_{\pi(\bvd n - 1)}$.
This permutation is well-defined because the proper subtrees of $t$
must be totally ordered under $\prec$:
if two such subtrees share the same maximum leaf label
(or indeed if they share any node or leaf at all),
they must be of different heights, or else they are the same subtree.

Conversely, if $\pi$ is a permutation on $0, \ldots, \bvd n - 1$,
we define the tree $T(\pi)$ as follows.
First, for $q$ in $Q$, we define
$\iota(q) = 1 + \max \{ \pi^{-1}(\bvd q + r) \;|\; r \in R \}$.
Since $\iota$ is injective,
we may unambiguously define, for $i$ in $0, \ldots, \bvd n$,
\[
 \tau_{i} =
  \begin{cases}
  \phi_{q}(\ldots, \tau_{\pi^{-1}(\bvd q + r)}, \ldots)
              & \parbox[t]{14em}{%
                if $i = \iota(q)$ for some $q \in Q$,} \\[0.5ex]
  \LEAF{n+i-\#\{ q \in Q \;|\; \iota(q) < i \}}
              & \parbox[t]{14em}{%
                otherwise.}
  \end{cases}
\]
We then take $T(\pi)$ to be $\tau_{\bvd n}$.

\exa{Let $d = 2$ and $n = 2$, and suppose the tree
$t = \PFXTREE{\PFXTREE{\LEAF{3}}{1}{\LEAF{2}}}{0}{\LEAF{4}}$ is given.
Its proper subtrees are
$\sigma_{0} = \PFXTREE{\LEAF{3}}{1}{\LEAF{2}}$,
$\sigma_{1} = \LEAF{4}$,
$\sigma_{2} = \LEAF{3}$, and
$\sigma_{3} = \LEAF{2}$.
These subtrees fall in the order
$\sigma_{3} \prec \sigma_{2} \prec \sigma_{0} \prec \sigma_{1}$,
which induces the permutation $\pi = \langle 3, 2, 0, 1 \rangle$;
thus we have $P(t) = \langle 3, 2, 0, 1 \rangle$.
Conversely, suppose the permutation $\pi = \langle 3, 2, 0, 1 \rangle$
is given.
The indices $\iota(q)$ are then
$\iota(0) = 1 + \max \{ \pi^{-1}(0), \pi^{-1}(1) \}
                = 1 + \max \{ 2, 3 \} = 4$, and
$\iota(1) = 1 + \max \{ \pi^{-1}(2), \pi^{-1}(3) \}
                = 1 + \max \{ 1, 0 \} = 2$.
We then have
$\tau_{0} = \LEAF{2}$;
$\tau_{1} = \LEAF{3}$;
$\tau_{2} = \PFXTREE{\tau_{\pi^{-1}(2)}}{1}%
                         {\tau_{\pi^{-1}(3)}} =
                 \PFXTREE{\tau_{1}}{1}{\tau_{0}}$;
$\tau_{3} = \LEAF{4}$; and
$\tau_{4} = \PFXTREE{\tau_{\pi^{-1}(0)}}{0}%
                         {\tau_{\pi^{-1}(1)}} =
                 \PFXTREE{\tau_{2}}{0}{\tau_{3}}$.
Hence $T(\pi) = \PFXTREE{\tau_{2}}{0}{\tau_{3}} =
 \PFXTREE{\PFXTREE{\LEAF{3}}{1}{\LEAF{2}}}{0}{\LEAF{4}}$.  $\square$
}

We now sketch a proof that $P$ and $T$ are inverses.
To begin, suppose we apply $T$ to a given $\pi$.
In the definition of $T$, each $\tau_{i}$ acquires a distinct root label,
hence all labels in $0, \ldots, \bvd n$ are represented.
Moreover, for $i$ in $0, \ldots, \bvd n - 1$,
the parent of $\tau_{i}$ is some $\tau_{i'}$
with $i < i' \le \bvd n$,
so by transitivity, each $\tau_{i}$ is a subtree of
$\tau_{\bvd n} = T(\pi)$.
By construction, we also have $\tau_{i} \prec \tau_{i+1}$ for each $i$,
hence $\tau_{0} \prec \ldots \prec \tau_{\bvd n}$.
Now suppose we apply $P$ to $T(\pi)$.
The $r$th child of node $q$ in $T(\pi)$ was defined to be
$\tau_{\pi^{-1}(\bvd q + r)}$,
so from the definition of $\sigma_{k}(t)$ it follows that
$\sigma_{\bvd q + r}(T(\pi)) = \tau_{\pi^{-1}(\bvd q + r)}$;
equivalently, we have $\sigma_{k}(T(\pi)) = \tau_{\pi^{-1}(k)}$
for $k$ in $0, \ldots, \bvd n - 1$.
Letting $i = \pi^{-1}(k)$, so that $\pi(i) = k$, we obtain
$\sigma_{\pi(i)}(T(\pi)) = \tau_{i}$ for $i$ in $0, \ldots, \bvd n - 1$,
which permits us to rewrite
$\tau_{0} \prec \ldots \prec \tau_{\bvd n - 1}$ as
$\sigma_{\pi(0)}(T(\pi)) \prec \ldots \prec \sigma_{\pi(\bvd n - 1)}(T(\pi))$.
Then by the definition of $P$, we have $P(T(\pi)) = \pi$.

Next consider $\hat{t} = T(P(t))$.
Let $\pi = P(t)$; then $P(\hat{t}) = P(T(P(t))) = P(T(\pi)) = \pi$ also.
Let us extend $\pi$ with $\pi(\bvd n) = \bvd n$,
and for all $t$ let $\sigma_{\bvd n}(t)$ denote $t$;
we then have
$\sigma_{\pi(0)}(t) \prec \ldots \prec \sigma_{\pi(\bvd n)}(t) = t$,
and similarly for $\hat{t}$.
We shall show by induction on $i$ that
$\sigma_{\pi(i)}(\hat{t}) = \sigma_{\pi(i)}(t)$,
and hence that $\hat{t} = t$.
If $\sigma_{\pi(i)}(\hat{t})$ is a node, it must have the form
$\phi_{q}(\ldots, \sigma_{dq+r}(\hat{t}), \ldots)$,
with $\sigma_{dq+r}(\hat{t}) \prec \sigma_{\pi(i)}(\hat{t})$
for each $r$ in $R$.
By
the inductive
hypothesis,
$\sigma_{dq+r}(\hat{t}) = \sigma_{dq+r}(t)$ for each $r$,
hence
$\phi_{q}(\ldots, \sigma_{dq+r}(\hat{t}), \ldots)$
is also a subtree of $t$.
Alternatively, if $\sigma_{\pi(i)}(\hat{t})$ is just a leaf $\LEAF{j}$,
this leaf must also occur as a subtree of $t$.
In either case, we deduce that
$\sigma_{\pi(i)}(\hat{t}) = \sigma_{\pi(i')}(t)$ for some $i' \ge i$.
Symmetrically, we have 
$\sigma_{\pi(i)}(t) = \sigma_{\pi(i'')}(\hat{t})$ for some $i'' \ge i$.
It follows that
$\sigma_{\pi(i)}(\hat{t}) \preceq \sigma_{\pi(i'')}(\hat{t}) =
 \sigma_{\pi(i)}(t) \preceq \sigma_{\pi(i')}(t) =
 \sigma_{\pi(i)}(\hat{t})$, and hence that
$\sigma_{\pi(i)}(\hat{t}) = \sigma_{\pi(i)}(t)$.

We remark without proof that $T$ and $P$ are exactly the mappings of
Sections~\ref{secperm2tree} and \ref{sectree2perm}.

\section{Conclusion}
\label{secconclude}

We have presented an algorithm for mapping
permutations to labeled trees,
as well as a variant of that algorithm
that performs the inverse mapping.
By establishing that these mappings are bijective,
we have shown that each of the factorials
in the formulas for the Catalan numbers and their $\bvd$-ary analogues
can be understood as a count of permutations.

\section*{Acknowledgments}

Discussions with
Peter Doyle helped me arrive at the present treatment
of this material.
Scott Daniels suggested notational improvements.

\end{document}